\documentclass[12pt,reqno]{article}
﻿
\usepackage[usenames]{color}
\usepackage[colorlinks=true,
linkcolor=webgreen, filecolor=webbrown,
citecolor=webgreen]{hyperref}
﻿
\definecolor{webgreen}{rgb}{0,.5,0}
\definecolor{webbrown}{rgb}{.6,0,0}
﻿
﻿
\usepackage{amssymb}
\usepackage{mathtools}
\usepackage{graphicx}
\usepackage{amscd}
\usepackage{lscape}
\usepackage{tikz}
\usepackage{tikz-cd}
\usepackage{pgfplots}
﻿
\usetikzlibrary{matrix}
\usetikzlibrary{fit,shapes}
\usetikzlibrary{positioning, calc}
\tikzset{circle node/.style = {circle,inner sep=1pt,draw, fill=white},
        X node/.style = {fill=white, inner sep=1pt},
        dot node/.style = {circle, draw, inner sep=5pt}
        }
\usepackage{tkz-fct}
﻿
\usepackage{amsthm}
\newtheorem{theorem}{Theorem}

\newtheorem{proposition}[theorem]{Proposition}
\newtheorem{corollary}[theorem]{Corollary}

﻿
\theoremstyle{definition}
\newtheorem{definition}[theorem]{Definition}
\newtheorem{example}[theorem]{Example}

﻿
\usepackage{float}
﻿
\usepackage{amsmath}
\usepackage{amsfonts}
\usepackage{latexsym}
\usepackage{epsf}
﻿
﻿

﻿
﻿
\newcommand{\seqnum}[1]{\href{http://oeis.org/#1}{\underline{#1}}}

\newcommand{\arxiv}[1]{\href{https://arxiv.org/abs/#1}{\texttt{arXiv:#1}}}
﻿
﻿
\begin{document}
﻿
﻿
\begin{center}
\vskip 1cm{\LARGE\bf Orthogonal polynomials on two or more intervals: the Sprugnoli context} \vskip 1cm \large
Paul Barry\\
School of Science\\
South East Technological University\\
Ireland\\
\href{mailto:pbarry@wit.ie}{\tt pbarry@wit.ie}
\end{center}
\vskip .2 in
﻿
﻿
﻿
\begin{abstract} We show how the group of Sprugnoli arrays on three elements can be used to define families of polynomials that are orthogonal on two intervals, and we indicate how the group can be extended in the hierarchy of Sprugnoli groups to yield families of orthogonal polynomials with more complex supports. In particular, we examine the defining three term recurrences and the Jacobi matrices of these families of orthogonal polynomials, while detailing the Sprugnoli matrix structure of the coefficient matrices of the polynomial families. We also consider integral representations of the principal moments on the defining intervals. Sprugnoli arrays are shown to relate to periodic recurrences where the parameters remain cyclically constant.
\end{abstract}
﻿
\section{Introduction}

The theory of orthogonal polynomials has a long and rich history, with many entry points, and many areas of application \cite{Chi, Gautschi, Geronimus, Ismail, Krushchev, Nikishin, Shohat, Sz}. In this paper, which looks principally at simple families of polynomials that are orthogonal on two or three intervals, we concentrate on four principal aspects. These are
\begin{itemize}
\item The structure of the coefficient array of the family
\item The corresponding three-term recurrence
\item The corresponding Jacobi matrix, and its role as the production matrix of the inverse of the coefficient array
\item The generating function of the  principal moment sequence of the family, and its integral representation on the intervals involved in the case of two intervals.
\end{itemize}

\section{Orthogonal polynomials supported on several intervals}
There is a long history of research into orthogonal polynomials defined on more than one interval \cite{Ac1, Ac2, Ac2, Geronimo, Peh1, Peh2}, and more generally on arcs and curves \cite{deB, Mil, Mil2, PehS1, PehS2}. The periodicity of the recurrence coefficients plays a central role in this context. This is evidenced in the associated Jacobi matrices. When the periodic coefficients remain constant (do not depend on $n$, for instance), major simplifications occur. This is the context of the present paper. Taking inspiration from the case of Chebyshev-type polynomials and their association with Riordan arrays \cite{BM, book2}, we examine orthogonal polynomials on two and three intervals defined by the extension of the Riordan group by the Sprugnoli groups \cite{Spru}. A major link between the elements of these groups and orthogonal polynomials is that these group elements are lower-triangular matrices determined by their production matrices, paralleling the fact that the Jacobi matrices of orthogonal polynomials are determined by the coefficient array of the polynomials.

\section{Recap on the Riordan group and Chebyshev-like orthogonal polynomials}
Given a family of orthogonal polynomials $P_n(x)$, the coefficient array $C$ of the family is the matrix  $C=(c_{n,k})$ where
$$c_{n,k}=[x^k] P_n(x).$$
Here, $[x^k]$ is the linear functional that acts on power series in $x$ to return the coefficient of $x^k$. As $P_n(x)$ is of degree $n$, the matrix $C$ is a lower-triangular matrix. Because $P_{n}(x)$ is of exact degree $n$, the matrix $C$ is invertible and the inverse of this matrix is associated with the moments of the family of orthogonal polynomials. In certain circumstances, the lower-triangular matrix $C$ can have a special structure which simplifies its analysis and applications. For instance, we can ask about the generating functions of the columns of the matrix $C$ (when viewed as sequences), and about those of the columns of its inverse. In certain cases, there is a structure to the sequence of these generating functions.

An example of this is given by looking at the Chebyshev polynomials of the second kind $U_n(x)$ \cite{Mason, Rivlin}. When we examine the coefficient array of the family of generalized scaled Chebyshev polynomials given by
$$P_n(x)=(\sqrt{b})^n U_n\left(\frac{x-a}{2\sqrt{b}}\right)-\alpha(\sqrt{b})^{n-1}U_{n-1}\left(\frac{x-a}{2\sqrt{b}}\right)-\beta(\sqrt{b})^{n-2}U_{n-2}\left(\frac{x-a}{2\sqrt{b}}\right)$$ we find that the generating functions of the columns of this matrix form a geometric series (of power series). This indicates that the coefficient array is an element of the Riordan group \cite{book1, book2, SGWW}. Specifically, it is given by the element
\begin{equation}\label{R1}
\left(\frac{1-\alpha x - \beta x^2}{1+ax+bx^2}, \frac{x}{1+ax+bx^2}\right)\end{equation} of the Riordan group \cite{book1, book2, SGWW}.

To describe this group,
we recall that $\mathcal{F}_0 = \{ \sum_{n=0}^{\infty}a_n x^n\,|\, a_0 \ne 0\}$ and that $\mathcal{F}_1 = \{ \sum_{n=0}^{\infty}a_n x^n\,|\, a_0 =0, a_1 \ne 0\}$ and in general
$$\mathcal{F}_r = \{\sum_{n=r}^{\infty}a_n x^n\,|\, a_r \ne 0\}.$$ For simplicity, we shall assume that the coefficients $a_n$ satisfy $a_n \in \mathbb{C}$. Elements of $\mathcal{F}_0$ are multiplicatively invertible, and elements of $\mathcal{F}_1$ are compositionally invertible.
﻿
The Riordan group is then the group of pairs $(g, f) \in \mathcal{F}_0 \times \mathcal{F}_1$ with the following product rule
\begin{equation}\label{product}(g(x), f(x))\cdot (u(x), v(x))=  (g(x)u(f(x)), v(f(x))\end{equation} and inverse
$$(g, f)^{-1}= \left(\frac{1}{g(\bar{f})}, \bar{f}(x)\right),$$ where $\bar{f}$ is the compositional inverse of $f \in \mathcal{F}_0$ (that is, $\bar{f}$ is the solution $u$ of the equation $f(u)=x$ for which $u(0)=0$). The identity of this group is $(1,x)$. To each element of this group we can associate in a unique way an invertible lower-triangular matrix $(t_{n,k})$ by means of
$$t_{n,k}=[x^n] g(x)f(x)^k,$$ where $[x^n]$ is the functional on $\mathbb{C}[[x]]$ that extracts the coefficient of $x^n$. Under this correspondence, the product (\ref{product}) corresponds to ordinary matrix multiplication. The columns of $(g,f)$ have their generating functions given by the geometric series of power series
$$g, gf, gf^2, gf^3, gf^4, \ldots.$$ The bi-variate generating function of the Riordan array $(g(x), f(x))$ is given by $\frac{g(x)}{1-yf(x)}$.
The product rule (\ref{product}) is the result of the following description of how a Riordan array acts on a power series.
\begin{equation}\label{Fun} (g(x), f(x))\cdot h(x)= g(x)h(f(x)). \end{equation} This weighted composition is conventionally referred to as the ``fundamental theorem of Riordan arrays''.

﻿If $M=(g,f)$ is a Riordan array, then the matrix $P_M=M^{-1}\overline{M}$, where $\overline{M}$ is the matrix $M$ with its first row removed, is called the \emph{production matrix} of the Riordan array $M$. We can recover $M=(t_{n,k})$ from $P_M$ by
$$t_{n,k}=\left(P_M^n\right)_{0,k}.$$
For instance, the production matrix of the inverse of the Riordan array (\ref{R1}) begins
$$\begin{pmatrix}
a + \alpha & 1 & 0 & 0 & 0 & 0 & 0 \\
b + \beta & a & 1 & 0 & 0 & 0 & 0 \\
0 & b & a & 1 & 0 & 0 & 0 \\
0 & 0 & b & a & 1 & 0 & 0 \\
0 & 0 & 0 & b & a & 1 & 0 \\
0 & 0 & 0 & 0 & b & a & 1 \\
0 & 0 & 0 & 0 & 0 & b & a
\end{pmatrix}.$$ This is the Jacobi matrix of the corresponding family of orthogonal polynomials. In this case, we call the inverse matrix $M=C^{-1}$ the \emph{moment matrix}, since the first column of this inverse matrix contains the moments of the orthogonal polynomials in question.
\begin{example} The Riordan array
$$C=\left(\frac{1}{1+2x+3x^2}, \frac{x}{1+2x+3x^2}\right)$$ begins
$$\begin{pmatrix}
1 & 0 & 0 & 0 & 0 & 0 & 0 \\
-2 & 1 & 0 & 0 & 0 & 0 & 0 \\
1 & -4 & 1 & 0 & 0 & 0 & 0 \\
4 & 6 & -6 & 1 & 0 & 0 & 0 \\
-11 & 4 & 15 & -8 & 1 & 0 & 0 \\
10 & -37 & -8 & 28 & -10 & 1 & 0 \\
13 & 72 & -66 & -40 & 45 & -12 & 1
\end{pmatrix}.$$
This is the coefficient array of the family of orthogonal polynomials $P_n(x)=(\sqrt{3})^n U_n\left(\frac{x-2}{2\sqrt{3}}\right)$. The inverse array $M$ (the moment array) begins
$$\begin{pmatrix}
1 & 0 & 0 & 0 & 0 & 0 & 0 \\
2 & 1 & 0 & 0 & 0 & 0 & 0 \\
7 & 4 & 1 & 0 & 0 & 0 & 0 \\
26 & 18 & 6 & 1 & 0 & 0 & 0 \\
106 & 80 & 33 & 8 & 1 & 0 & 0 \\
452 & 365 & 170 & 52 & 10 & 1 & 0 \\
1999 & 1692 & 861 & 304 & 75 & 12 & 1
\end{pmatrix}.$$
Here, the sequence \seqnum{A122871} $$1,2,7,26,106,452,1999,\ldots,$$ is the moment sequence of the family of orthogonal polynomials. Note that where known, we refer to sequences by their entry number in the On-Line Encyclopedia of Integer Sequences (OEIS) \cite{SL1, SL2}.

The production matrix $P_M$ of this inverse array begins
$$\begin{pmatrix}
2 & 1 & 0 & 0 & 0 & 0 & 0 \\
3 & 2 & 1 & 0 & 0 & 0 & 0 \\
0 & 3 & 2 & 1 & 0 & 0 & 0 \\
0 & 0 & 3 & 2 & 1 & 0 & 0 \\
0 & 0 & 0 & 3 & 2 & 1 & 0 \\
0 & 0 & 0 & 0 & 3 & 2 & 1 \\
0 & 0 & 0 & 0 & 0 & 3 & 2
\end{pmatrix}.$$ Thus we have $P_{-1}(x)=0$, $P_0(x)=1$, and
$$P_n(x)=(x-2)P_{n-1}(x)-3P_{n-2}.$$
The elements of $P_M$ define both the inverse matrix $M$ and the original coefficient array $C$. For instance, the element $33$ in $M$ is given by  $$33=1\cdot 18+2\cdot 6+3\cdot 1.$$ In general, we have
$$m_{n,k}=1 \cdot m_{n-1,k-1} + 2 \cdot m_{n-1,k}+ 3\cdot m_{n-1,k+1}.$$
This illustrates a defining characteristic of all Riordan arrays: each internal element (an element not in the first column) is a linear combination of the element in the row above it, starting by going one up and one to the left. The same linear combination coefficients hold for all elements. This is a consequence of the linear recurrence characterization of Riordan arrays \cite{book2}.

With regard to the original coefficient array $C$, we have, for instance, that the element $15$ visible above is given by
$$15 = 1\cdot 6 + (-2)\cdot (-6) + (-3)\cdot 1.$$
Note that here, we use $(1,-2,-3)$ as multipliers. This is also a consequence of a linear recurrence characterization of Riordan arrays, but can be seen directly by
\begin{align*}c_{n,k}&=[x^k]P_n(x)\\&=[x^k](x-a)P_{n-1}(x)-b P_{n-2}(x)\\
&=[x^k]xP_{n-1}(x)-a[x^k]P_{n-1}(x)-b[x^k]P_{n-2}(x)\\
&=[x^{k-1}]P_{n-1}(x)-a[x^k]P_{n-1}(x)-b[x^k]P_{n-2}(x)\\
&=c_{n-1,k-1}-a\cdot c_{n-1,k}-b\cdot c_{n-2,k}.\end{align*}
We have
$$\left(\frac{1}{1+2x+3x^2}, \frac{x}{1+2x+3x^2}\right)^{-1}=\left(\frac{1-2x-\sqrt{1-4x-8x^2}}{6x^2},\frac{1-2x-\sqrt{1-4x-8x^2}}{6x}\right).$$
Thus $\frac{1-2x-\sqrt{1-4x-8x^2}}{6x^2}$ is the generating function of the moment sequence. Using the Stieltjes-Perron transform \cite{Chi}, this gives us the integral representation of the moments
$$\mu_n=m_{n,0}=\frac{1}{\pi} \int_{2-2\sqrt{3}}^{2+2\sqrt{3}} {{\frac{x^n \sqrt{-x^2+4x+8}}{6}dx}}.$$
The support of the corresponding measure is then the interval
$$[2-2\sqrt{3},2+2\sqrt{3}].$$
We note further that the generating function $\mu(x)$ of the moments $\mu_n$ can be represented by the Jacobi continued fraction \cite{Wall}
$$\cfrac{1}{1-2x-\cfrac{3x^2}{1-2x-\cfrac{3x^2}{1-2x-\cdots}}}.$$ Such sequences play an important role in the enumeration of lattice paths \cite{Flajolet, Viennot}. For instance, the sequence $1,2,7,26,\ldots$ counts Motzkin paths where the flat steps have $2$ possible colors and the up steps have $3$ possible colors.
\end{example}
In general, we have
$$\left(\frac{1-\alpha x - \beta x^2}{1+ax+bx^2}, \frac{x}{1+ax+bx^2}\right)^{-1}=$$
\begin{scriptsize}
$$\left(\frac{(b+\beta)\sqrt{1-2ax+(a^2-4b)x^2}+(a(b-\beta)+2b\alpha)-b+\beta}{2(x^2(a^2\beta+a\alpha(\beta-b)-b^2-b(\alpha^2+2\beta)-\beta^2)+x(\alpha(b-\beta)-2\alpha \beta)+\beta)},\frac{1-ax-\sqrt{1-2ax-(a^2-4b)x^2}}{2bx}\right).$$
\end{scriptsize} Again, the Stieltjes-Perron transform leads us to the weight function
$$\frac{(b+\beta)\sqrt{-x^2+2ax-a^2+4b}}{2(\beta x^2+x(\alpha(b-\beta)-2\alpha \beta)+a^2\beta+a \alpha(\beta-b)-b^2-b(\alpha^2+2\beta)-\beta^2}$$ on the interval
$$[a-2\sqrt{b},a+2\sqrt{b}].$$
For completeness we mention that the coefficient array of the Chebyshev polynomials of the first kind $T_n(x)$ is not a Riordan array, but a so-called ``almost'' Riordan array, an element of the group of almost Riordan arrays (a supergroup of the Riordan group) \cite{almost}.

\section{Orthogonal polynomials on two intervals: initial Sprugnoli examples}

In order to motivate our discussion of the Sprugnoli group in the context of orthogonal polynomials on two intervals, we begin with two examples. We will give more specific information on the Sprugnoli group, and its elements, Sprugnoli arrays, in the next section.
\begin{example} We consider the family of orthogonal polynomials $P_n(x)$ defined by $P_{-1}(x)=0$, $P_0(x)=1$, and
\begin{align*} P_n(x)&=(x-2)P_{n-1}(x)-3 P_{n-2}(x),\quad \text{if $n$ is even}\\
P_n(x)&=(x+2)P_{n-1}(x)-5 P_{n-2}(x), \quad \text{if $n$ is odd}.\end{align*}
The coefficient matrix $C$ of this family of orthogonal polynomials begins
$$\begin{pmatrix}
1 & 0 & 0 & 0 & 0 & 0 & 0 \\
-2 & 1 & 0 & 0 & 0 & 0 & 0 \\
-9 & 0 & 1 & 0 & 0 & 0 & 0 \\
24 & -12 & -2 & 1 & 0 & 0 & 0 \\
93 & 0 & -21 & 0 & 1 & 0 & 0 \\
-258 & 129 & 48 & -24 & -2 & 1 & 0 \\
-981 & 0 & 330 & 0 & -33 & 0 & 1
\end{pmatrix}.$$
The general term $c_{n,k}$ of this matrix can be shown to be given by the following expression.
\begin{equation}\label{Spuxn}c_{n,k} = [x^n] g(x) f_1(x)^{k \bmod 2} (x f_2(x))^{\lfloor \frac{k}{2} \rfloor}.\end{equation}
Here, the triple of power series $(g(x), f_1(x), f_2(x))$ is given by
$$(g(x),f_1(x), f_2(x))=\left(\frac{1-2x+3x^2}{1+12x^2+15x^4}, \frac{x}{1-2x+3x^2}, \frac{x}{1+12x^2+15x^4}\right).$$
This triple, and its matrix representation, is an example of a \emph{Sprugnoli array}.

The inverse $M$ of this array begins
$$\begin{pmatrix}
1 & 0 & 0 & 0 & 0 & 0 & 0 \\
2 & 1 & 0 & 0 & 0 & 0 & 0 \\
9 & 0 & 1 & 0 & 0 & 0 & 0 \\
18 & 12 & 2 & 1 & 0 & 0 & 0 \\
96 & 0 & 21 & 0 & 1 & 0 & 0 \\
192 & 159 & 42 & 24 & 2 & 1 & 0 \\
1179 & 0 & 363 & 0 & 33 & 0 & 1
\end{pmatrix}.$$ We regard the first column sequence $1,2,9,18,96,\ldots$ as the principal moments of this family of orthogonal polynomials.

The production matrix $P$ of the inverse matrix $M$ begins
$$P=\begin{pmatrix}
2 & 1 & 0 & 0 & 0 & 0 & 0 \\
5 & -2 & 1 & 0 & 0 & 0 & 0 \\
0 & 3 & 2 & 1 & 0 & 0 & 0 \\
0 & 0 & 5 & -2 & 1 & 0 & 0 \\
0 & 0 & 0 & 3 & 2 & 1 & 0 \\
0 & 0 & 0 & 0 & 5 & -2 & 1 \\
0 & 0 & 0 & 0 & 0 & 3 & 2
\end{pmatrix}.
$$
This is the periodic Jacobi matrix of the family of $P_n(x)$ of orthogonal polynomials. Being periodic of order two, we expect the family of orthogonal polynomials to be supported on two intervals.

From this Jacobi matrix, we obtain that the moment sequence $1,2,9,18,96,\ldots$ has its generating function given by the continued fraction
$$\mu(x)=\cfrac{1}{1-2x-\cfrac{5x^2}{1+2x-\cfrac{3x^2}{1-2x-\cfrac{5x^2}{1+2x-\cfrac{3x^2}{1-2x-\cdots}}}}}.$$
To obtain a closed-form expression for this generating function, we solve the equation
$$u=\cfrac{1}{1-2x-\cfrac{5x^2}{1+2x-3x^2u}}.$$
We find the following generating function for the principal moments $\mu_n$
$$\mu(x)=\frac{1-6x^2-\sqrt{1-24x^2+84x^4}}{6x^2(1-2x)}.$$
Using the Stieltjes-Perron transform again, we find the following integral representation for the moments $\mu_n$.
$$\mu_n = \frac{1}{\pi}\left\{ \int_{\sqrt{12-\sqrt{60}}}^{\sqrt{12+\sqrt{60}}} x^n \frac{\sqrt{-x^4+24x^2-84}}{6(x-2)}\,dx- \int_{-\sqrt{12+\sqrt{60}}}^{-\sqrt{12-\sqrt{60}}} x^n \frac{\sqrt{-x^4+24x^2-84}}{6(x-2)}\,dx\right\}.$$
Thus we see that the family of orthogonal polynomials $P_n(x)$ is supported on two intervals, namely
$$\left[-\sqrt{12+\sqrt{60}}, -\sqrt{12-\sqrt{60}}\right] \quad \text{and}\quad \left[\sqrt{12-\sqrt{60}},\sqrt{12+\sqrt{60}}\right].$$
\end{example}
\begin{example} In this example, we consider the family of orthogonal polynomials $Q_n(x)$ defined by $Q_{-1}(x)=0$, $Q_0(x)=1$, and
\begin{align*} Q_n(x)&=(x-2)Q_{n-1}(x)-5 Q_{n-2}(x),\quad \text{if $n$ is odd}\\
Q_n(x)&=(x+2)Q_{n-1}(x)-3 Q_{n-2}(x), \quad \text{if $n$ is even}.\end{align*}
The coefficient array $C_Q$ of this family of orthogonal polynomials is then given by the Sprugnoli array
$$\left(\frac{1-2x+5x^2}{1+12x^2+15x^4}, \frac{x}{1-2x+5x^2}, \frac{x}{1+12x^2+15x^4}\right).$$
We find that the principal moments $\mu_n$ begin $1,2,7,14,64,128,\ldots$, with generating function given by the Jacobi continued fraction
$$\mu(x)=\cfrac{1}{1-2x-\cfrac{3x^2}{1+2x-\cfrac{5x^2}{1-2x-\cfrac{3x^2}{1+2x-\cdots}}}}.$$
We deduce that the sequence of moments $\mu_n$ has the integral representation
$$\mu_n= \frac{1}{\pi}\left\{ \int_{\sqrt{12-\sqrt{60}}}^{\sqrt{12+\sqrt{60}}} x^n \frac{\sqrt{-x^4+24x^2-84}}{10(x-2)}\,dx- \int_{-\sqrt{12+\sqrt{60}}}^{-\sqrt{12-\sqrt{60}}} x^n \frac{\sqrt{-x^4+24x^2-84}}{10(x-2)}\,dx\right\}+\frac{2}{5}\cdot 2^n.$$
Thus in this case, there is an additional atomic weight at $x=2$.
\end{example}

\section{The Sprugnoli group}

In this section, we briefly summarize the details of the Sprugnoli group that we need. The Sprugnoli group \cite{Spru} is a generalization of the Riordan group \cite{book1, book2, SGWW}. Elements of the Sprugnoli group are defined by a triple of power series $(g, f_1, f_2)$, while elements of the Riordan group are defined by a pair of power series $(g, f)$.

We have seen that the product law in the Riordan group follows from the following result, called the ``fundamental theorem of Riordan arrays'', which details how a Riordan array operates on a power series. We have
$$(g(x), f(x))\cdot h(x)= g(x)h(f(x)).$$
﻿
A similar law holds for so-called ``vertically stretched'' Riordan arrays \cite{LI, book2}. These are matrices corresponding to pairs of power series of the form $(g, xf(x))$ where as usual, $f(x) \in \mathcal{F}_1$, and so $xf(x) \in \mathcal{F}_2$. Again, we have a weighted composition action on power series given by
$$ (g(x), xf(x))\cdot h(x)= g(x)h(xf(x)).$$
The generating function of the stretched Riordan array $(g, xf(x))$  is given by
$$(g(x),xf(x))\cdot \frac{1}{1-xy}=\frac{g(x)}{1-yxf(x)}.$$ In particular, the generating function of the row sums of the stretched Riordan array $(g,xf)$ will have generating function $\frac{g(x)}{1-xf(x)}$.
\begin{example} The stretched Riordan array $\left(\frac{1}{1-x}, \frac{x^2}{1-x-x^2}\right)$ begins
$$\left(\begin{array}{ccccccccc}
1 & 0 & 0 & 0 & 0 & 0 & 0 & 0 & 0 \\
1 & 0 & 0 & 0 & 0 & 0 & 0 & 0 & 0 \\
1 & 1 & 0 & 0 & 0 & 0 & 0 & 0 & 0 \\
1 & 2 & 0 & 0 & 0 & 0 & 0 & 0 & 0 \\
1 & 4 & 1 & 0 & 0 & 0 & 0 & 0 & 0 \\
1 & 7 & 3 & 0 & 0 & 0 & 0 & 0 & 0 \\
1 & 12 & 8 & 1 & 0 & 0 & 0 & 0 & 0 \\
1 & 20 & 18 & 4 & 0 & 0 & 0 & 0 & 0 \\
1 & 33 & 38 & 13 & 1 & 0 & 0 & 0 & 0
\end{array}\right).$$ Applying this matrix to the positive Fibonacci numbers \seqnum{A000045}, we obtain
$$\left(\begin{array}{ccccccccc}
1 & 0 & 0 & 0 & 0 & 0 & 0 & 0 & 0 \\
1 & 0 & 0 & 0 & 0 & 0 & 0 & 0 & 0 \\
1 & 1 & 0 & 0 & 0 & 0 & 0 & 0 & 0 \\
1 & 2 & 0 & 0 & 0 & 0 & 0 & 0 & 0 \\
1 & 4 & 1 & 0 & 0 & 0 & 0 & 0 & 0 \\
1 & 7 & 3 & 0 & 0 & 0 & 0 & 0 & 0 \\
1 & 12 & 8 & 1 & 0 & 0 & 0 & 0 & 0 \\
1 & 20 & 18 & 4 & 0 & 0 & 0 & 0 & 0 \\
1 & 33 & 38 & 13 & 1 & 0 & 0 & 0 & 0
\end{array}\right)\left( \begin{array}{c}1\\1\\2\\3\\5\\8\\13\\21\\34\end{array}\right)=\left( \begin{array}{c}1\\1\\2\\3\\7\\14\\32\\69\\159\end{array}\right),$$ where the last sequence has generating function
$$\left(\frac{1}{1-x},\frac{x^2}{1-x-x^2}\right)\cdot \frac{1}{1-x-x^2}=\frac{\frac{1}{1-x}}{1-\frac{x^2}{1-x-x^2}-(\frac{x^2}{1-x-x^2})^2},$$ which is equal to
$$\frac{(1-x-x^2)^2}{(1-x)(1-2x+3x^3+x^4)}.$$
\end{example}

\begin{definition} The \emph{Sprugnoli group} is the set with elements $(g, f_1, f_2)$ with $g(x) \in \mathcal{F}_0$, $f_1(x) \in \mathcal{F}_1$, and $f_2(x) \in \mathcal{F}_1$ and
$f_2 \in xR[[x^2]]$ (thus $f_2$ is an odd power series). The element $(g, f_1, f_2)$ of this group has the matrix representation $S=(s_{n,k})$ where
$$s_{n,k} = [x^n] g(x)f_1(x)^{k \bmod 2} (xf_2(x))^{\lfloor \frac{k}{2} \rfloor}.$$
\end{definition}
Thus we have
$$
s_{n,k} =
\begin{cases}
  [x^n]g(x)(xf_2(x))^m, & k=2m, \\
  [x^n]g(x)f_1(x)(xf_2(x))^m,  & k=2m+1.
\end{cases}
$$
The columns of the matrix $S$ have generating functions given by the following products of generating functions
\begin{align*}g(x), g(x)f_1(x),&\, g(x)(xf_2(x)), g(x)f_1(x)(xf_2(x)), g(x)(xf_2(x))^2,\\
& g(x)f_1(x)(xf_2(x))^2, g(x)(xf_2(x))^3, g(x)f_1(x)(xf_2(x))^3,\ldots\end{align*} We can represent this sequence by the following schema.
$$
\begin{array}{cccccccc}
g &g & g & g & g & g & g & g\\
1 &f_1 & 1 & f_1 & 1 & f_1 & 1 & f_1 \\
1 &1 & f_2 & f_2 & f_2^2 & f_2^2 & f_2^3 & f_2^3 \\
1 & 1 & x & x & x^2 & x^2 & x^3 & x^3 \\
- & - & - & - & - & - & - & -\\
0 & 1 & 2 & 3 & 4 & 5 & 6 & 7 \\
\end{array}
$$ Here, the final row gives the order of the product of the terms above it, showing that this associated matrix is lower-triangular.

Thus
\begin{itemize}
\item $g$ is present in all columns
\item $f_1$ only appears in every other column, always to the power of one
\item $xf_2$ appears to the power of $\lfloor \frac{k}{2} \rfloor$: $0,0,1,1,2,2,\ldots$

\begin{align*}
g&:1,1,1,1,1,1,1\ldots\\
f_1&:0,1,0,1,0,1,0,1,0,1\ldots\\
xf_2&:0,0,1,1,2,2,3,3,\ldots \end{align*}
\end{itemize}

We can represent the matrix $S$ as the sum of two matrices as follows. The first matrix is the matrix whose columns are generated by the power series
$$g(x), 0, g(x)(xf_2(x)), 0, g(x)(xf_2(x))^2, \ldots,$$ which is a horizontal aeration of the stretched Riordan array $(g(x), xf_2(x))$, and the second matrix is the matrix whose columns are generated by the power series
$$0,g(x)f_1(x),0, g(x)f_1(x)(xf_2(x)), 0, g(x)f_1(x)(xf_2(x))^2,0,\ldots,$$ which is the horizontal aeration of the stretched Riordan array $(g(x)f_1(x), xf_2(x))$ (but note that this matrix starts with a zero row, due to the fact that $f_1 \in \mathcal{F}_1$).\
﻿
\begin{example} We consider the Sprugnoli array $\left(\frac{1}{1-x}, \frac{x(1+x)}{1-x}, \frac{x}{1-x^2}\right)$. This array begins
$$\left(\begin{array}{ccccccccc}
1 & 0 & 0 & 0 & 0 & 0 & 0 & 0 & 0 \\
1 & 1 & 0 & 0 & 0 & 0 & 0 & 0 & 0 \\
1 & 3 & 1 & 0 & 0 & 0 & 0 & 0 & 0 \\
1 & 5 & 1 & 1 & 0 & 0 & 0 & 0 & 0 \\
1 & 7 & 2 & 3 & 1 & 0 & 0 & 0 & 0 \\
1 & 9 & 2 & 6 & 1 & 1 & 0 & 0 & 0 \\
1 & 11 & 3 & 10 & 3 & 3 & 1 & 0 & 0 \\
1 & 13 & 3 & 15 & 3 & 7 & 1 & 1 & 0 \\
1 & 15 & 4 & 21 & 6 & 13 & 4 & 3 & 1
\end{array}\right).$$ This is the sum of the matrices
$$\left(\begin{array}{ccccccccc}
1 & 0 & 0 & 0 & 0 & 0 & 0 & 0 & 0 \\
1 & 0 & 0 & 0 & 0 & 0 & 0 & 0 & 0 \\
1 & 0 & 1 & 0 & 0 & 0 & 0 & 0 & 0 \\
1 & 0 & 1 & 0 & 0 & 0 & 0 & 0 & 0 \\
1 & 0 & 2 & 0 & 1 & 0 & 0 & 0 & 0 \\
1 & 0 & 2 & 0 & 1 & 0 & 0 & 0 & 0 \\
1 & 0 & 3 & 0 & 3 & 0 & 1 & 0 & 0 \\
1 & 0 & 3 & 0 & 3 & 0 & 1 & 0 & 0 \\
1 & 0 & 4 & 0 & 6 & 0 & 4 & 0 & 1
\end{array}\right)+
\left(\begin{array}{ccccccccc}
0 & 0 & 0 & 0 & 0 & 0 & 0 & 0 & 0 \\
0 & 1 & 0 & 0 & 0 & 0 & 0 & 0 & 0 \\
0 & 3 & 0 & 0 & 0 & 0 & 0 & 0 & 0 \\
0 & 5 & 0 & 1 & 0 & 0 & 0 & 0 & 0 \\
0 & 7 & 0 & 3 & 0 & 0 & 0 & 0 & 0 \\
0 & 9 & 0 & 6 & 0 & 1 & 0 & 0 & 0 \\
0 & 11 & 0 & 10 & 0 & 3 & 0 & 0 & 0 \\
0 & 13 & 0 & 15 & 0 & 7 & 0 & 1 & 0 \\
0 & 15 & 0 & 21 & 0 & 13 & 0 & 3 & 0
\end{array}\right).$$
The first is a horizontal aeration of the stretching of the  Riordan array $\left(\frac{1}{1-x},\frac{x}{1-x^2}\right)$, while the second one
is a horizontal aeration of the stretching of the Riordan array $\left(\frac{1}{1-x} \frac{x(1+x)}{1-x},\frac{x}{1-x^2}\right)=\left(\frac{x(1+x)}{(1-x)^2},\frac{x}{1-x^2}\right)$ with an extra initial column of zeros.
\end{example}

A special element of this set is given by $(1,x,x)$. Its columns are generated by the sequence $1,x,x^2,x^3,\ldots$ and so its matrix representation is given by the usual identity matrix.
﻿
Elements of this set are referred to as \emph{Sprugnoli arrays}, in reference to their matrix representation.
\begin{proposition}
The bi-variate generating function of the Sprugnoli array $(g, f_1, f_2)$ is given by
\begin{equation}\label{bi}\frac{g(x)}{1-y^2xf_2(x)}+\frac{yg(x)f_1(x)}{1-y^2xf_2(x)}=\frac{g(x)(1+yf_1(x))}{1-y^2xf_2(x)}.\end{equation}
\end{proposition}
\begin{proof} The array is given by the sum of the horizontal aeration of the stretched Riordan array $(g(x), xf_2(x))$ and the horizontal aeration of the stretched Riordan array $(g(x)f_1(x), xf_2(x))$. Interpreting this in terms of generating functions gives the following results.
\end{proof}
\begin{corollary} The row sums of the Sprugnoli array $(g, f_1, f_2)$ have generating function $\frac{g(x)(1+f_1(x))}{1-xf_2(x)}$.
\end{corollary}
\begin{proof} This results by setting $y=1$ in the generating function of the array.
\end{proof}
\begin{corollary} The diagonal sums of the Sprugnoli array $(g, f_1, f_2)$ have generating function $\frac{g(x)(1+xf_1(x))}{1-x^3f_2(x)}$.
\end{corollary}
\begin{proof} This results by setting $y=x$ in the generating function of the array.
\end{proof}
﻿
\section{The ``fundamental theorem'' of Sprugnoli arrays}

﻿What is conventionally called the ``fundamental theorem'' of Riordan arrays is the prescription of how a Riordan array operates on a power series. Similarly, we have the ``fundamental theorem'' of Sprugnoli arrays, which details how a Sprugnoli array operates on a power series. This operation uses the even and odd bisections of the power series being operated on. We have the following result.
\begin{proposition} Let $(g, f_1, f_2)$ be a Sprugnoli array, and let $h(x)$ be a power series $h_n(x)=\sum_{n=0}^{\infty} a_n x^n$. We let
$$h^e(x)=\sum_{n=0}^{\infty} a_{2n}x^n \quad\text{and}\quad h^o(x)=\sum_{n=0}^{\infty} a_{2n+1}x^n$$ be the even and odd bisections of $h(x)$. Then we have
$$(g(x), f_1(x), f_2(x))\cdot h(x)=g(x)h^e(xf_2(x))+g(x)f_1(x)h^o(xf_2(x)).$$
\end{proposition}
\begin{proof} We partition the matrix representing $(g, f_1, f_2)$ as
$$(g,0,g(xf_2),0,g(xf_2)^2,0,g(xf_2)^3,0,\ldots) + (0, gf_1, 0, gf_1(xf_2), 0, gf_1(xf_2)0,\ldots).$$
These matrices operate, respectively, on $(a_0,0,a_2,0,a_4,0,\ldots)$ and $(0,a_1,0,a_3,0,a_5,0,\ldots)$.
In terms of Riordan arrays, this is equivalent, after compression, to the stretched Riordan array $(g, xf_2)$ operating on $h^e$, and the
stretched Riordan array $(gf_1, xf_2)$ operating on $h^o(x)$. The result follows from this.
\end{proof}
We recall that $g^e(x)=\frac{g(\sqrt{x})+g(-\sqrt{x})}{2}$ and $g^o(x)=\frac{g(\sqrt{x})-g(-\sqrt{x})}{2\sqrt{x}}$.
We can write this in matrix form as
$$\frac{1}{2}
\begin{pmatrix}
1 & 0 \\
0 & x^{-1/2}
\end{pmatrix}
\begin{pmatrix}
1 & 1 \\
1 & -1
\end{pmatrix}
\begin{pmatrix}
g\!\left(x^{1/2}\right) \\
g\!\left(-x^{1/2}\right)
\end{pmatrix}
=
\begin{pmatrix}
g^e(x) \\
g^o(x)
\end{pmatrix}.$$
\begin{example} We have
\begin{align*}(1, x, x) \cdot h(x)&= 1.h^e(x^2)+1.x.h^o(x^2)\\
&=h^e(x^2)+xh^o(x^2)\\
&=h(x).\end{align*}
\end{example}
\section{The product rule}
We are now in a position to say what we mean by the product of two Sprugnoli arrays. We have the following.
\begin{definition} We define the product of two Sprugnoli arrays $(g, f_1, f_2)$ and $(u, v_1, v_2)$ as follows.
\begin{equation}\label{prod}(g, f_1, f_2) \cdot (u, v_1, v_2)= \left((g,f_1,f_2)\cdot u, \frac{(g,f_1,f_2)\cdot u v_1}{(g,f_1,f_2)\cdot u}, \frac{1}{x}.\frac{(g,f_1,f_2)\cdot uxv_2}{ (g,f_1,f_2)\cdot u}\right).\end{equation}
\end{definition}
This can be simplified to the following.
\begin{equation}\label{prod2}(g, f_1, f_2) \cdot (u, v_1, v_2)= \left((g,f_1,f_2)\cdot u, \frac{(g,f_1,f_2)\cdot u v_1}{(g,f_1,f_2)\cdot u}, \frac{1}{x} \sqrt{xf_2}v_2(\sqrt{xf_2})\right).\end{equation}
Furthermore, the product rule (\ref{prod}), when we look to the matrix representation of Sprugnoli arrays, coincides with matrix multiplication.
﻿
\section{The inverse of a Sprugnoli array}
﻿
We have
$$(g,f_1,f_2)^{-1}=\left(1,\frac{x-f_1^e(\left(\overline{\sqrt{xf_2}}\right)^2)}{f_1^o(\left(\overline{\sqrt{xf_2}}\right)^2)},\frac{1}{x}\left(\overline{\sqrt{xf_2}}\right)^2\right)\cdot \left(\frac{1}{g},x,x\right).$$ We define generating functions $r_1(x)$ and $r_2(x)$ by
$$(1, r_1(x), r_2(x))=\left(1,\frac{x-f_1^e(\left(\overline{\sqrt{xf_2}}\right)^2)}{f_1^o(\left(\overline{\sqrt{xf_2}}\right)^2)},\frac{1}{x}\left(\overline{\sqrt{xf_2}}\right)^2\right).$$

\section{The production matrix of a Sprugnoli array}
A defining feature of Riordan arrays is that they have a sequence characterization. For instance, in a Riordan array $M=(g,f)$, any matrix element $t_{n,k}$ not in the first column will obey a linear recurrence equation involving the elements in the row above it, starting with the element $t_{n-1,k-1}$ above and to the left of it. The same recurrence coefficients are valid for each such element, irrespective of the row it is in. This means that the production matrix $P$ of a Riordan array has a special simple structure, where $P=M^{-1}\overline{M}$, where $\overline{M}$ is the matrix $M$ with its top row removed \cite{Prod}. The elements in the first column of $M$ normally follow a distinct recurrence. This results in two sequences of recurrence coefficients, the $Z$ sequence, for the first column elements, and the $A$ sequence, for all the other elements \cite{SEQ}. The production array will then have a structure as in the following matrix, where the obvious pattern continues.
$$\left(\begin{array}{rrrrrrr}
z_0 & a_0 & 0 & 0 & 0 &  0 & \cdots \\
z_1 & a_1 & a_0 & 0 & 0 & 0 & \cdots \\
z_2 & a_2 & a_1 &a_0 & 0 & 0 & \cdots \\
z_3 & a_3 & a_2 & a_1 & a_0 & 0 & \cdots \\
z_4 & a_4 & a_3 & a_2 & a_1 & a_0& \cdots \\
z_5 & a_5 & a_4 & a_3 & a_2 & a_1 & \cdots \\
\vdots & \vdots & \vdots & \vdots & \vdots & \vdots & \ddots
\end{array}\right).$$
For Sprugnoli arrays, a similar situation arises, but this time, we have two $A$-sequences (we designate them the $A$ sequence and the $B$-sequence), one for elements in even columns, and one for elements in odd columns. Thus the form of the production array for a Sprugnoli array is given by
\begin{equation}\label{Stripe}P=\left(\begin{array}{rrrrrrr}
z_0 & a_0 & 0 & 0 & 0 &  0 & \cdots \\
z_1 & a_1 & b_0 & 0 & 0 & 0 & \cdots \\
z_2 & a_2 & b_1 &a_0 & 0 & 0 & \cdots \\
z_3 & a_3 & b_2 & a_1 & b_0 & 0 & \cdots \\
z_4 & a_4 & b_3 & a_2 & b_1 & a_0& \cdots \\
z_5 & a_5 & b_4 & a_3 & b_2 & a_1 & \cdots \\
\vdots & \vdots & \vdots & \vdots & \vdots & \vdots & \ddots
\end{array}\right),\end{equation} giving it a ``striped'' appearance.
We let $A(x)$ be the generating function of the second column of the production matrix (beginning $a_0, a_1, a_2,\ldots$). If the third column of the production matrix begins $0,b_0, b_1, b_2,\ldots$, then we let $B(x)$ be the generating function of the sequence $b_0, b_1, b_2,\ldots$. We have the following result.
\begin{proposition}
We have
\begin{align*} A(x)&=(1,r_1,r_2)\cdot \frac{f_1(x)}{x},\\
B(x)&=\frac{1}{x} (1,r_1,r_2)\cdot f_2(x).\end{align*}
\end{proposition}
We also have the following necessary condition for a production matrix to be that of a Sprugnoli array.
\begin{proposition}
Let $A(x)$ and $xB(x)$ correspond to the generating functions of the second and third columns of the production array of the Sprugnoli array $(g, f_1, f_2)$. Then $A(x)+B(x)$ is an even power series.
\end{proposition}
﻿
\section{Orthogonal polynomials on two intervals: the Sprugnoli case}
We are now in a position to exploit the Sprugnoli group formulation to present a general result concerning ``constant coefficient'' orthogonal polynomials on two intervals. By the designation  ``constant coefficient'' we mean that the periodic coefficients of the three-term recurrence remain constant. The rows and columns of the corresponding Jacobi matrix will then repeat two by two. We have the following result.
\begin{proposition} The coefficient array of the  family of orthogonal polynomials $P_n(x)$ defined by $P_{-1}(x)=0$, $P_0(x)=1$, and
\begin{align*} P_n(x)&=(x-a)P_{n-1}(x)-c P_{n-2}(x),\quad \text{if $n$ is odd}\\
P_n(x)&=(x+a)P_{n-1}(x)-b P_{n-2}(x), \quad \text{if $n$ is even},\end{align*} is given by the Sprugnoli array
$$(g(x), f_1(x), f_2(x))=\left(\frac{1-ax+cx^2}{1+(a^2+b+c)x^2+bcx^4}, \frac{x}{1-ax+cx^2}, \frac{x}{1+(a^2+b+c)x^2+bcx^4}\right),$$
with $f_2(x)=g(x)f_1(x)$.
The inverse of the coefficient array (the moment array) is given by $(w(x), s_1(x), s_2(x))$ where
\begin{align*}
\mu(x)&=\frac{1+ax}{1-(a^2+b-c)x^2}c\left(\frac{cx^2(1-a^2x^2)}{(1-(a^2+b-c)x^2)^2}\right),\\
s_1(x)&=\frac{x(1-ax)}{1-(a^2-b+c)x^2}c\left(\frac{bx^2(1-a^2x^2)}{(1-(a^2-b+c)x^2)^2}\right),\\
﻿s_2(x)&=\mu(x)s_1(x),\end{align*}
where $$c(x)=\frac{1-\sqrt{1-4x}}{2x}$$ is the generating function of the Catalan numbers \seqnum{A000108}.
The production matrix of the moment matrix (the Jacobi matrix of the family $P_n(x)$) begins
\begin{equation}\label{Jac0} \begin{pmatrix}
a &  1 &  0 &  0 &  0 &  0 & 0 \\
b & -a &  1 &  0 &  0 &  0 & 0 \\
0 &  c &  a &  1 &  0 &  0 & 0 \\
0 &  0 &  b & -a &  1 &  0 & 0 \\
0 &  0 &  0 &  c &  a &  1 & 0 \\
0 &  0 &  0 &  0 &  b & -a & 1 \\
0 &  0 &  0 &  0 &  0 &  c & a
\end{pmatrix}. \end{equation}
The principal moments $\mu_n$ have generating function $w(x)$ and weight function
$$\frac{\sqrt{-x^4+2x^2(a^2+b+c)-a^4-2a^2(b+c)-(b-c)^2}}{2(x-a)}$$ supported on the two intervals
$$\left[-\sqrt{a^2+b+c+2\sqrt{bc}}, -\sqrt{a^2+b+c-2\sqrt{bc}}\right], \left[\sqrt{a^2+b+c-2\sqrt{bc}}, \sqrt{a^2+b+c+2\sqrt{bc}}\right].$$
\end{proposition}
\begin{proof}
The theory of Sprugnoli arrays tells us that the production matrix of the array $(g(x), f_1(x), f_2(x))$ above is given by the Jacobi matrix (\ref{Jac0}). The theory of orthogonal polynomials then tells us that this Jacobi matrix corresponds to the family $P_n(x)$ of orthogonal polynomials. The generating function of the principal moments $\mu_n$ is given by the generating function of the first column of the inverse array. An application of the Stieltjes-Perron transform to this generating function then yields the weight function for the moments $\mu_n$ on the two supporting intervals.
\end{proof}
\begin{corollary} The principal moments $\mu_n$ of the family of orthogonal polynomials $P_n(x)$ have their generating function $\mu(x)$ given by the Jacobi continued fraction
$$\mu(x)=\frac{1}{1-ax-\cfrac{bx^2}{1+ax-\cfrac{cx^2}{1-ax-\cfrac{bx^2}{1+ax-\cdots}}}}.$$
\end{corollary}
\begin{corollary}
The Hankel transform $h_n=|\mu_{i+j}|_{0 \le i,j \le n}$ of the moment sequence $\mu_n$ is given by
$$h_n = b^{\lfloor \frac{(n+1)^2}{4} \rfloor} c^{\lfloor \frac{n^2}{4} \rfloor}.$$
\end{corollary}
We also have that $s_1(x)$ has its generating function given by
$$s_1(x)=\frac{x}{1+ax-\cfrac{cx^2}{1-ax-\cfrac{bx^2}{1+ax-\cfrac{cx^2}{1-ax-\cdots}}}}.$$
We deduce that the Hankel transform of the expansion of $s_1(x)/x$ is given by $c^{\lfloor \frac{(n+1)^2}{4} \rfloor} b^{\lfloor \frac{n^2}{4} \rfloor}$.

\begin{corollary} The generating function of the sequence of orthogonal polynomials $P_n(x)$ is given by 
$$\frac{1+(x-a)y+cy^2}{1-(x^2-(a^2+b+c))y^2+bcy^4}.$$
\end{corollary}
\begin{proof} We use the formula  for the bi-variate generating function of a Sprugnoli array (\ref{bi}) to conclude that 
$$P_n(x)=[y^n] \frac{1+(x-a)y+cy^2}{1-(x^2-(a^2+b+c))y^2+bcy^4}.$$
\end{proof}

\begin{example} We take the case $a=1$, $b=1$ and $c=2$. Thus we consider the family of orthogonal polynomials $P_n(x)$ with $P_{-1}(x)=0$, $P_0(x)=1$, and
\begin{align*} P_n(x)&=(x-1)P_{n-1}(x)-2 P_{n-2}(x),\quad \text{if $n$ is odd}\\
P_n(x)&=(x+1)P_{n-1}(x)-P_{n-2}(x), \quad \text{if $n$ is even}.\end{align*} The coefficient array of this family is given by the Sprugnoli array
$$(g(x), f_1(x), f_2(x))=\left(\frac{1-x+2x^2}{1+4x^2+2x^4}, \frac{x}{1-x+2x^2}, \frac{x}{1+4x^2+2x^4}\right).$$
We find that
$$\mu(x)=(1+x)c(2x^2(1-x^2))=\frac{1-\sqrt{1-8x^2+8x^4}}{4x^2(1-x)}$$ from which we deduce that
$$\mu_n=\frac{1}{\pi}\left\{\int_{4-\sqrt{8}}^{4+\sqrt{8}} x^n \frac{\sqrt{-x^4+8x^2-8}}{4(x-1)}\,dx-\int_{-(4+\sqrt{8})}^{-(4-\sqrt{8})} x^n \frac{\sqrt{-x^4+8x^2-8}}{4(x-1)}\,dx \right\}+\frac{1}{2}.$$
Here, $\frac{1}{2}=\frac{1}{2}(1)^n$ is an atomic weight at $x=1$ occasioned by the factor $x-1$ in the denominator. The sequence $\mu_n$ begins $1,1,2,2,6,6,24,24,112,\ldots$. This is a doubling of the sequence \seqnum{A118376}. A continued fraction expression for $\mu(x)$ is given by
$$\mu(x)=\cfrac{1}{1-x-\cfrac{x^2}{1+x-\cfrac{2x^2}{1-x-\cfrac{x^2}{1+x-\cfrac{2x^2}{1-x-\cdots}}}}}.$$
As a doubled sequence, we also have
$$\mu(x)=\cfrac{1+x}{1-x^2-\cfrac{x^2}{1-\cfrac{2x^2}{1-x^2-\cfrac{x^2}{1-\cfrac{2x^2}{1-x^2-\cfrac{x^2}{1-\cdots}}}}}}.$$
\end{example}
\begin{example} To contrast with the above example, we look at the case $a=1$, $b=2$ and $c=1$. Thus we are considering the family of orthogonal polynomials $P_n(x)$ with $P_{-1}(x)=0$, $P_0(x)=1$, and
\begin{align*} P_n(x)&=(x-1)P_{n-1}(x)-P_{n-2}(x),\quad \text{if $n$ is odd}\\
P_n(x)&=(x+1)P_{n-1}(x)-2P_{n-2}(x), \quad \text{if $n$ is even}.\end{align*} The coefficient of this family is given by the Sprugnoli array
$$(g(x), f_1(x), f_2(x))=\left(\frac{1-x+x^2}{1+4x^2+2x^4}, \frac{x}{1-x+x^2}, \frac{x}{1+4x^2+2x^4}\right).$$
In this case, we find that
$$\mu(x)=\frac{1+x}{1-2x^2}c\left(\frac{x^2(1-x^2)}{(1-2x^2)^2}\right)=\frac{1-2x^2-\sqrt{1-8x^2+8x^4}}{2x^2(1-x)}.$$ We deduce that
$$\mu_n=\frac{1}{\pi}\left\{\int_{4-\sqrt{8}}^{4+\sqrt{8}} x^n \frac{\sqrt{-x^4+8x^2-8}}{2(x-1)}\,dx-\int_{-(4+\sqrt{8})}^{-(4-\sqrt{8})} x^n \frac{\sqrt{-x^4+8x^2-8}}{2(x-1)}\,dx \right\}.$$ The sequence $\mu_n$ in this case begins $1,1,3,3,11,11,47,47,223,\ldots$. It is a doubling of the sequence \seqnum{A174347}. A continued fraction expression for $\mu(x)$ is given by
$$\mu(x)=\cfrac{1}{1-x-\cfrac{2x^2}{1+x-\cfrac{x^2}{1-x-\cfrac{2x^2}{1+x-\cfrac{x^2}{1-2x-\cdots}}}}}.$$
As a doubled sequence, we also have
$$\mu(x)=\cfrac{1+x}{1-x^2-\cfrac{2x^2}{1-\cfrac{x^2}{1-x^2-\cfrac{2x^2}{1-\cfrac{x^2}{1-x^2-\cfrac{2x^2}{1-\cdots}}}}}}.$$
\end{example}
We note that the sequences $1,2,6,24,\ldots$ and $1,3,11,47,\ldots$ count coloured Schroeder paths.

In general, we have that
$$\mu(x)=(1+ax)g(x^2),$$ where
\begin{align*}g(x)&=\frac{1}{1-(a^2+b-c)x}c\left(\frac{cx(1-a^2x)}{(1-(a^2+b-c)x)^2}\right)\\
&=\cfrac{1}{1-a^2x-\cfrac{bx}{1-\cfrac{cx}{1-a^2x-\cfrac{bx}{1-\cfrac{cx}{1-a^2x-\cdots}}}}}.\end{align*}
This Thron continued fraction \cite{Thron} is the generating function for colored Schroeder paths with no level paths at odd level, level paths of color (weight) $a^2$ at even levels, and alternating up steps of color $b$ and $c$. In addition, we have that $g(x)$ has the following Jacobi continued fraction expression
$$g(x)=\cfrac{1}{1-(a^2+b)x-\cfrac{bcx^2}{1-(a^2+b+c)x-\cfrac{bcx^2}{1-(a^2+b+c)x-\cfrac{bcx^2}{1-(a^2+b+c)x-\cdots}}}},$$
which enumerates colored Motzkin paths \cite{Flajolet, Viennot}.

\section{The next Sprugnoli group and orthogonal polynomials}
It is possible to generalize the Sprugnoli group to the group that is made up of elements $(g,f_1,\ldots,f_n)$ for general $n$, where $g(x) \in \mathcal{F}_0$, $f_1,\ldots,f_{n-1}\in\mathcal{F}_1$, and $f_n \in x\mathbb{C}[[x^n]]$.  In this section, we shall apply the group of elements $(g, f_1, f_2, f_3)$ to the study of certain families of orthogonal polynomials defined on three sets.
For this, we specify that $f_3(x) \in x \mathbf{F}[[x^3]]$ \cite{Spru}, and we consider the group of matrices $C=(c_{n,k})$ where
$$c_{n,k}= [x^n]g(x)f_1(x)^{\frac{2}{3}\left(1-\cos\left(\frac{2 \pi k}{3}\right)\right)}f_2(x)^{\frac{2}{3}\left(\frac{1}{2}+\cos\left(\frac{2 \pi (k+1)}{3}\right)\right)}(x^2f_3(x))^{\lfloor \frac{k}{3} \rfloor}.$$ Such arrays follow the following schema, which indicates the powers of the constituent elements $g(x), f_1(x), f_2(x), f_3(x)$ in the generating functions of the successive columns of the array.
\begin{align*}
g&:1,1,1,1,1,1,1\ldots\\
f_1&:0,1,1,0,1,1,0,1,1,0,\ldots\\
f_2&:0,0,1,0,0,1,0,0,1,0,0,1,\ldots\\
x^2f_3&:0,0,0,1,1,1,2,2,2,3,\ldots \end{align*}

To define the ``fundamental theorem of the $3$-Sprugnoli group'', we need to calculate the three trisections of a power series. For a generating function $h(x)$, we can calculate its trisections $h^{(1)}(x)$, $h^{(2)}(x)$ and $h^{(3)}{(x)}$ using the $3 \times 3$ discrete Fourier transform \cite{DFT} as follows. We let $\omega=\exp(2\pi i/3)$. Then
$$\frac{1}{3}
\begin{pmatrix}
1 & 0 & 0\\
0 & x^{-1/3} & 0\\
0 & 0 & x^{-2/3}
\end{pmatrix}
\begin{pmatrix}
1 & 1 & 1 \\
1 & \omega^2 & \omega\\
1 & \omega & \omega^2
\end{pmatrix}
\begin{pmatrix}
h\!\left(x^{1/3}\right) \\
h\!\left(\omega x^{1/3}\right)\\
h\!\left(\omega^2 x^{1/3}\right)
\end{pmatrix}
=
\begin{pmatrix}
h^{(1)}(x) \\
h^{(2)}(x) \\
h^{(3)}(x)
\end{pmatrix}.$$
The fundamental theorem for this group, describing how an element $(g, f_1, f_2, f_3)$ acts on a power series $h(x)$, is given by
$$(g,f_1,f_2,f_3)\cdot h(x)=g(x)h^{(1)}(x^2f_3(x))+g(x)f_1(x)h^{(2)}(x^2f_3(x))+g(x)f_1(x)f_2(x)h^{(3)}(x^2f_3(x)).$$
Note that we can write this in the symbolic form
$$(g,f_1,f_2,f_3)\cdot h(x)=\begin{pmatrix} (g,x^2f_3), & (gf_1,x^2f_3), & (gf_1f_2,x^2f_3) \end{pmatrix} \begin{pmatrix} h^{(1)}\\h^{(2)}\\h^{(3)}\end{pmatrix}.$$

We now wish to associate matrices of the type $(g(x), f_1(x), f_2(x), f_3(x))$ with the coefficient arrays of orthogonal polynomials that are defined on three intervals.
\begin{example}
We consider that Sprugnoli matrix given by
\begin{scriptsize}
$$(g(x),f_1(x),f_2(x),f_3(x))=\left(\frac{1-x-x^2+7x^3-14x^4}{1+7x^3-28x^6},\frac{x(1-2x-7x^2)}{1-x-x^2+7x^3-14x^4},\frac{x}{1-2x-7x^2},\frac{x}{1+7x^3-28x^6}\right).$$
\end{scriptsize}
The production matrix of the inverse $(\mu(x),s_1(x),s_2(x),s_3(x))$ of this array begins
$$\begin{pmatrix}
1 & 1 & 0 & 0 & 0 & 0 & 0 & 0 \\
2 & 1 & 1 & 0 & 0 & 0 & 0 & 0 \\
0 & 2 & -2 & 1 & 0 & 0 & 0 & 0 \\
0 & 0 & -7 & 1 & 1 & 0 & 0 & 0 \\
0 & 0 & 0 & 2 & 1 & 1 & 0 & 0 \\
0 & 0 & 0 & 0 & 2 & -2 & 1 & 0 \\
0 & 0 & 0 & 0 & 0 & -7 & 1 & 1 \\
0 & 0 & 0 & 0 & 0 & 0 & 2 & 1
\end{pmatrix}.$$
The first column of the inverse array begins
$$1, 1, 3, 7, 21, 49, 119, 315, 749,\ldots,$$ with generating function
$$\mu(x)=\cfrac{1}{1-x-\cfrac{2x^2}{1-x-\cfrac{2x^2}{1+2x+\cfrac{7x^2}{1-x-\cfrac{2x^2}{1-x-\cfrac{2x^2}{1+2x+\cfrac{7x^2}{1-x-\cdots}}}}}}}.$$
Equivalently, we have
$$\mu(x)=\frac{\sqrt{1-14x^3+161x^6}-7x^3+14x^2-1}{14x^2(1-2x-x^2)}.$$
We have that
$$(g(x),f_1(x), f_2(x), f_3(x))^{-1}=(\mu(x),s_1(x),s_2(x),s_3(x))$$ where
\begin{align*}
s_1(x)&=\cfrac{x}{1-x-\cfrac{2x^2}{1+2x+\cfrac{7x^2}{1-x-\cfrac{2x^2}{1-x-\cfrac{2x^2}{1+x+\cdots}}}}}\\
&=\frac{1+4x^2+x^3-\sqrt{1-14x^3+161x^6}}{4x(1+x-4x^2)},\end{align*}
\begin{align*}
s_2(x)&=\cfrac{x}{1+x+\cfrac{7x^2}{1-x-\cfrac{2x^2}{1-x-\cfrac{2x^2}{1+x+\cdots}}}}\\
&=\frac{1+4x^2-11x^3-\sqrt{1-14x^3+161x^6}}{4x(1+x+5x^2)},\end{align*}
and
$$s_3(x)=\frac{\mu(x)s_1(x)s_2(x)}{x}.$$
The matrix $C$ is then the coefficient array for the family of orthogonal polynomials $P_n(x)$ for which $P_{-1}(x)=0$, $P_0(x)=1$, $P_1(x)=x-1$, and
\begin{align*}
P_n(x)&=(x+2)P_{n-1}(x)-2P_{n-2}(x) \quad\text{if $n \bmod 3=0$},\\
P_n(x)&=(x-1)P_{n-1}(x)+7 P_{n-2}(x) \quad\text{if $n \bmod 3=1$},\\
P_n(x)&=(x-1)P_{n-1}(x)-2P_{n-2}(x) \quad\text{if $n \bmod 3=2$}.\end{align*}
The three contours of integration in this case lie in the complex plane, with endpoints given by the six zeros of the polynomial  $x^6-14x^3+161$, which occur in conjugate pairs on the circle centred at $(0,0)$ and of radius $161^{\frac{1}{6}}$. Specifically, these zeros are
$(7-4i\sqrt{7})^{1/3}, -(7+4i\sqrt{7})^{1/3}, (-1)^{2/3}(7-4i\sqrt{7})^{1/3}, (7+4i\sqrt{7})^{1/3},-(-1)^{1/3}(7+4i\sqrt{7})^{1/3}, (-1)^{2/3}(7+4i\sqrt{7})^{1/3}$.
\begin{figure}[htbp] \centering \includegraphics[width=0.4\textwidth]{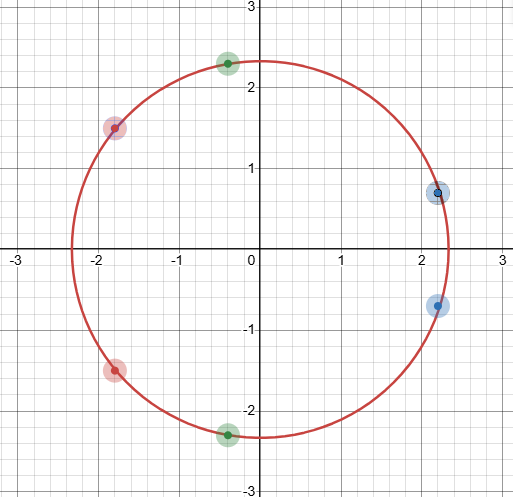} \caption{Zeros of $x^6 - 14x^3 + 161$} \label{fig:zeros} \end{figure}
This suggests that we are in the presence of polynomials orthogonal on several arcs \cite{PehS1}.
\end{example}
We have the following proposition, describing when a Sprugnoli array $(g(x), f_1(x), f_2(x), f_3(x))$ is the coefficient array of a family of orthogonal polynomials.
\begin{proposition} The Sprugnoli array $(g(x), f_1(x), f_2(x), f_3(x))$ where
\begin{align*}
g(x)&=\frac{1-ux-(a-uv)x^2+v(a+b+u^2+uv+v^2)x^3-a(a+b+u^2+uv+v^2)x^4}{1+(a(u+2v)+b(v-u)+v^3)x^3-ab(a+b+u^2+uv+v^2)x^6},\\
f_1(x)&=\frac{x(1-(u+v)x-(a+b+u^2+uv+v^2)x^2)}{1-ux-(a-uv)x^2+v(a+b+u^2+uv+v^2)x^3-a(a+b+u^2+uv+v^2)x^4},\\
f_2(x)&=\frac{x}{1-(u+v)x-(a+b+u^2+uv+v^2)x^2},\\
f_3(x)&=\frac{g(x)f_1(x)f_2(x)}{x},\end{align*} is the coefficient array of the family $P_n(x)$ of orthogonal polynomials where $P_{-1}(x)=0$, $P_0(x)=1$, $P_1(x)=x-u$, and
\begin{align*}
P_n(x)&=(x+u+v)P_{n-1}(x)-bP_{n-2}(x) &\text{if $n \bmod 3=0$},\\
P_n(x)&=(x-u)P_{n-1}(x)-c P_{n-2}(x) &\text{if $n \bmod 3=1$},\\
P_n(x)&=(x-v)P_{n-1}(x)-aP_{n-2}(x) &\text{if $n \bmod 3=2$},\end{align*}
where $$c=-(a+b+u^2+uv+v^2).$$
\end{proposition}
\begin{proof}
The production matrix of the inverse matrix begins (with a continuing pattern)
\begin{equation}\label{Jac}\begin{pmatrix}
u & 1 & 0 & 0 & 0 & 0 & 0 & 0 \\
a & v & 1 & 0 & 0 & 0 & 0 & 0 \\
0 & b & -(u+v) & 1 & 0 & 0 & 0 & 0 \\
0 & 0 & c & u & 1 & 0 & 0 & 0 \\
0 & 0 & 0 & a & v & 1 & 0 & 0 \\
0 & 0 & 0 & 0 & b & -(u+v) & 1 & 0 \\
0 & 0 & 0 & 0 & 0 &  c & u & 1 \\
0 & 0 & 0 & 0 & 0 &  0 & a & v \end{pmatrix}.\end{equation}
This is the Jacobi matrix of the family of orthogonal polynomials $P_n(x)$. The first column sequence of the inverse matrix (the principal moment sequence) has its generating function given by
$$w(x)=\cfrac{1}{1-ux-\cfrac{ax^2}{1-vx-\cfrac{bx^2}{1+(u+v)x-\cfrac{cx^2}{1-ux-\cfrac{ax^2}{1-vx-\cdots}}}}},$$ from which we deduce that
the family of orthogonal polynomials $P_n(x)$ is supported on the (complex) intervals determined by the zeros of the polynomial
\begin{scriptsize}
$$x^6 - 2x^3(a(u + 2v) + b(v - u) + v^3) + a^2(4b + u^2 + 4uv + 4v^2) + 2a(2b^2 + b(u^2 + uv + 4v^2) + v^3(u + 2v)) + b^2(u^2 - 2uv + v^2) + 2bv^3(v - u) + v^6.$$
\end{scriptsize}
The condition $c=-(a+b+u^2+uv+v^2)$ is to ensure that the matrix $C$ is indeed a Sprugnoli array.
\end{proof}
\begin{example} The condition $c=-(a+b+u^2+uv+v^2)$ is necessary for the coefficient array to be a Sprugnoli array. We give an example of a family of orthogonal polynomials $P_n(x)$ whose Jacobi matrix is of the above form (\ref{Jac}) but where we do not have $c=-(a+b+u^2+uv+v^2)$.
We let $P_{-1}(x)=0$, $P_0(x)=1$, $P_1(x)=x-1$, and
\begin{align*}
P_n(x)&=(x+2)P_{n-1}(x)-2P_{n-2}(x) &\text{if $n \bmod 3=0$},\\
P_n(x)&=(x-1)P_{n-1}(x)-P_{n-2}(x) &\text{if $n \bmod 3=1$},\\
P_n(x)&=(x-1)P_{n-1}(x)-2P_{n-2}(x) &\text{if $n \bmod 3=2$}.\end{align*}
The Jacobi matrix for this family of orthogonal polynomials begins
$$\begin{pmatrix}
1 & 1 & 0 & 0 & 0 & 0 & 0 & 0 \\
2 & 1 & 1 & 0 & 0 & 0 & 0 & 0 \\
0 & 2 & -2 & 1 & 0 & 0 & 0 & 0 \\
0 & 0 & 1 & 1 & 1 & 0 & 0 & 0 \\
0 & 0 & 0 & 2 & 1 & 1 & 0 & 0 \\
0 & 0 & 0 & 0 & 2 & -2 & 1 & 0 \\
0 & 0 & 0 & 0 & 0 & 1 & 1 & 1 \\
0 & 0 & 0 & 0 & 0 & 0 & 2 & 1
\end{pmatrix},$$ with a repeating pattern. The principal moments $\mu_n$ of this family of orthogonal polynomials is the sequence that begins
$$1, 1, 3, 7, 21, 49, 151, 347, 1101, 2485, 8139,\ldots,$$ whose generating function is given by
\begin{align*}w(x)&=\cfrac{1}{1-x-\cfrac{2x^2}{1-x-\cfrac{2x^2}{1+2x-\cfrac{x^2}{1-x-\cfrac{2x^2}{1-x-\cdots}}}}}\\
&=\frac{1-6x^2-x^3-\sqrt{1-16x^2+2x^3+64x^4-16x^5-15x^6}}{2x^2(1-2x-x^2)}.
\end{align*}
The support of this family is then determined by the location of the zeros of the polynomial $x^6 - 16x^4 + 2x^3 + 64x^2 - 16x - 15$, which are all real. We have (in order)
\begin{align*}
x_1&=\frac{\sqrt{5}}{2}-\frac{3}{2},\\
x_2&=\frac{-\sqrt{5}}{2}-\frac{3}{2},\\
x_3&=-\frac{4\sqrt{6}}{3}\sin\left(\frac{1}{3}\tan^{-1} \left(\frac{15\sqrt{4119}}{1373}\right)+\frac{\pi}{3}\right),\\
x_4&=\frac{4\sqrt{6}}{3}\cos\left(\frac{1}{3}\cot^{-1} \left(\frac{15\sqrt{4119}}{1373}\right)\right),\\
x_5&=\frac{4\sqrt{6}}{3}\sin\left(\frac{1}{3}\tan^{-1} \left(\frac{15\sqrt{4119}}{1373}\right)\right),\\
x_6&=3.
\end{align*}
Then we have the integral representation
$$\mu_n=\frac{1}{\pi}\left\{\int_{x_1}^{x_2} x^n W(x)dx-\int_{x_3}^{x_4} x^n W(x)dx+\int_{x_5}^{x_6} x^n W(x)dx\right\},$$ where
$$W(x)=\frac{\sqrt{-x^6+16x^4-2x^3-64x^2+16x+15}}{2(x^2-2x-1)}.$$ Thus this family of orthogonal polynomials is supported on three intervals (but it does not arise from a Sprugnoli array).
\end{example}

\section{Conclusions} We have approached the study of orthogonal polynomials on two intervals by asking what type of structure does this confer on the coefficient array of the corresponding family of orthogonal polynomials. We have shown that in the so-called ``constant coefficient'' case, this structure is given by a Sprugnoli array $(g,f_1,f_2)$. That is, it is an element of the Sprugnoli group (of order one). This allows us to bring to bear all the theory of the Sprugnoli group. We have also initiated a study of the connections between orthogonal polynomials and elements $(g,f_1,f_2,f_3)$ of the Sprugnoli group of order two. This case appears to be more challenging, and awaits further study.

\bigskip
\hrule
\bigskip
\noindent 2020 {\it Mathematics Subject Classification}:
Primary 15B36; Secondary 05A15, 11B83, 15A30, 20H25.
﻿
\noindent \emph{Keywords:} Orthogonal polynomials, Riordan group, Sprugnoli group, generating function, production matrix, Jacobi matrix.
﻿
﻿
﻿

﻿
﻿
\end{document}